\documentclass[11pt,reqno]{amsart}

\usepackage{amssymb}
\usepackage{amsthm}
\usepackage{amsmath,amsfonts,amssymb}
\usepackage{graphics,color}
\usepackage{enumerate}
\usepackage{float}

\newtheorem{theorem}{Theorem}[section]
\newtheorem{lemma}[theorem]{Lemma}
\newtheorem{proposition}[theorem]{Proposition}

\newtheorem{definition}[theorem]{Definition\rm}
\newtheorem{remark}{\it Remark\/}

\newcommand*{\R}{\ensuremath{\mathbb{R}}}
\renewcommand*{\S}{\ensuremath{\mathcal{S}}}

\newcommand*{\N}{\ensuremath{\mathbb{N}}}

\newcommand*{\supp}{\ensuremath{\mathrm{supp\,}}}
\newcommand*{\dist}{\ensuremath{\mathrm{dist\,}}}
\renewcommand*{\div}{\ensuremath{\mathrm{div\,}}}

\newcommand*{\e}{\ensuremath{\varepsilon}}
\newcommand{\eps}{\varepsilon}
\def\rn#1{\mathbb{R}^{#1}}

\begin{document}

\title{The Euler equations as a differential inclusion}

\author{Camillo De Lellis}
\address{Institut f\"ur Mathematik, Universit\"at Z\"urich, CH-8057 Z\"urich}
\email{camillo.delellis@math.unizh.ch}

\author{L\'aszl\'o Sz\'ekelyhidi Jr.}
\address{Departement Mathematik, ETH Z\"urich, CH-8092 Z\"urich}
\email{szekelyh@math.ethz.ch}

\begin{abstract}
In this paper we propose a new point of view on weak solutions of the Euler equations, describing the motion of an ideal incompressible fluid in $\R^n$ with $n\geq 2$. We give a reformulation of the Euler equations as a differential inclusion, and in this way we obtain transparent proofs of several celebrated results of V.~Scheffer and A.~Shnirelman concerning the non-uniqueness of weak solutions and the existence of energy--decreasing solutions. Our results are stronger because they work in any dimension and yield bounded velocity and pressure. 
\end{abstract}

\maketitle

\section{Introduction}

Consider the Euler equations in $n$ space dimensions, describing the motion of an ideal incompressible fluid, 
\begin{equation}\label{euler}
\begin{split}
\partial_tv+\div(v\otimes v)+\nabla p - f&= 0\\
\div v&=0\, .
\end{split}
\end{equation}
Classical (i.e.~sufficiently smooth) solutions of the Cauchy problem 
exist locally in time for sufficiently regular
initial data and driving forces (see Chapter 3.2 in \cite{MajdaBook}). In two dimensions such existence results
are available also for global solutions (e.g.~Chapters 3.3 and 8.2 in 
\cite{MajdaBook} and the references therein).
Classical solutions of
Euler's equations with $f=0$ conserve the energy, that is
$t\mapsto \int |v (x,t)|^2\, dx$ is a constant function.
Hence the energy space for \eqref{euler} is
$L^\infty_t (L^2_x)$. 

A recurrent issue in the modern theory of PDEs
is that one needs to go beyond classical solutions,
in particular down to the energy space
(see for instance \cite{DafermosBook,DipernaMajda1,MajdaBook,TaoBook}).
A divergence--free vector field $v\in L^2_{loc}$ is a 
{\em weak solution} of 
\eqref{euler} if
\begin{equation}\label{distrib}
\int \big(v \partial_t \varphi + \langle v\otimes v, 
\nabla \varphi\rangle + \varphi\cdot f\big)\, dx\, dt \;=\; 0\\
\end{equation}
for every test function $\varphi\in C^\infty_c 
(\R^n_x\times \R_t, \R^n)$ with ${\rm div}\, \varphi = 0$.
It is well--known that then the pressure is determined up to
a function depending only on time (see \cite{TemamBook}).
In the case of Euler
strong motivation for considering weak solutions comes also from
mathematical physics, especially the theory of turbulence laid down by Kolmogorov
in 1941 \cite{ChorinBook,FrischBook}. A celebrated criterion of Onsager related to Kolmogorov's theory says, roughly speaking, that dissipative weak solutions cannot have a H\"older exponent greater than $1/3$ 
(see \cite{ConstantinETiti,DuchonRobert,Eyink,Onsager}). 
It is therefore of interest to construct weak solutions with limited regularity. 

Weak solutions are not unique. In a well--known paper \cite{Scheffer93} 
Scheffer constructed a surprising example of a weak solution to \eqref{euler} with 
compact support in space and time when $f=0$ and $n=2$. Scheffer's
proof is very long and complicated and a simpler
construction was later given by Shnirelman in 
\cite{Shnirelman1}. 
However, Shnirelman's proof is still quite difficult. In this paper we obtain
a short and elementary proof of the following
theorem. 
\begin{theorem}\label{main}
Let $f=0$. There exists $v\in L^{\infty}(\R^n_x\times\R_t;\R^n)$ 
and $p\in L^{\infty}(\R^n_x\times\R_t)$
solving \eqref{euler}
in the sense of distributions, such that $v$ 
is not identically zero, and $\supp v$ and $\supp p$ 
are compact in space-time $\R^n_x\times\R_t$.
\end{theorem}

In mathematical physics weak solutions to the Euler equations that dissipate energy underlie the Kolmogorov theory of turbulence. In another groundbreaking paper \cite{Shnirelmandecrease} Shnirelman proved the existence of $L^2$ distributional solutions with $f=0$ and energy which decreases in time.
His methods are completely unrelated to those in \cite{Scheffer93} and \cite{Shnirelman1}. In contrast, the following extension of his existence theorem is  a simple corollary of our construction.
\begin{theorem}\label{t:decrease}
There exists $(v,p)$ as in Theorem \ref{main} such that, in addition:
\begin{itemize}
\item $\int |v (x,t)|^2\, dx =1$ for almost every $t\in ]-1,1[$,
\item $v (x,t)=0$ for $|t|>1$.
\end{itemize}
\end{theorem}
Our method has several interesting features.
First of all, our approach  fits nicely in the well--known framework of L.~Tartar
for the analysis of oscillations in linear
partial differential systems coupled with nonlinear pointwise constraints \cite{Diperna85,KMS02,Tartar79,Tartar82}.
Roughly speaking, Tartar's framework amounts to a plane--wave analysis localized in physical space, 
in contrast with Shnirelman's method in \cite{Shnirelman1}, 
which is based rather on a wave analysis in Fourier space.
In combination with Gromov's 
convex integration or with Baire category arguments, Tartar's approach leads to a well
understood mechanism for generating irregular oscillatory solutions to 
differential inclusions (see \cite{Bernd,KMS02,MS99}). 

Secondly, the velocity field we construct belongs to the energy space
 $L^\infty_t (L^2_x)$. This was not the case
for the solutions in \cite{Scheffer93,Shnirelman1}, and it was a natural question whether weak solutions in the energy space were unique.
Our first theorem shows that even
higher summability assumptions
of $v$ do not rule out such pathologies. 
The pressure in \cite{Scheffer93,Shnirelman1} is only a distribution solving
\eqref{euler}. In our construction
$p$ is actually the potential--theoretic
solution of
\begin{equation}\label{eqpress}
- \Delta p \;=\; \partial^2_{x_ix_j} (v^iv^j) - 
\partial_{x_i} f_i\, .
\end{equation}
However, being bounded,
it has slightly better regularity than the 
$BMO$ given by the classical
estimates for \eqref{eqpress}. 

Next, our point of view reveals connections between the apparently unrelated 
constructions of Scheffer and Shnirelman.
Shnirelman considers sequences of 
driving forces $f_k$ converging to $0$
in some negative Sobolev space. 
In particular he shows that for a suitable choice of $f_k$
the corresponding solutions of \eqref{euler} 
converge in $L^2$ to a nonzero solution of \eqref{euler} with $f=0$.
Scheffer builds his solution
by iterating a certain piecewise
constant construction at small scales.
On the one hand both our proof and Scheffer's proof are based on oscillations localized
in physical space. On the other hand, our proof
gives as an easy byproduct the following approximation 
result in Shnirelman's spirit.

\begin{theorem}\label{main2} All the solutions $(v,p)$ constructed in the
proofs of  Theorem \ref{main} and in Theorem \ref{t:decrease} have 
the following property. There exist three sequences
$\{v_k\}, \{f_k\}, \{p_k\} \subset C^\infty_c$ 
solving \eqref{euler} such that
\begin{itemize}
\item $f_k$ converges to $0$ in $H^{-1}$,
\item $\|v_k\|_\infty+ \|p_k\|_\infty$ 
is uniformly bounded,
\item $(v_k,p_k)\to (v,p)$
in $L^q$ for every $q<\infty$.
\end{itemize}
\end{theorem}

Our results give interesting information on which kind of additional
(entropy) condition could restore uniqueness of solutions. As already remarked,  belonging to
the energy space is not sufficient.
In fact, in view of our method of construction, there is strong evidence that neither energy--decreasing nor energy--preserving solutions are unique.
In a forthcoming paper we plan to investigate this issue, and also the class of initial data for which our method yields energy--decreasing solutions. 

The rest of the paper is organized as follows. In Section \ref{Tartar} we carry out the plane wave analysis of the Euler equations in the spirit of Tartar, and we formulate the core of our construction (Proposition \ref{p:lam}). In Section \ref{s:planewaves} we prove Proposition \ref{p:lam}. In Section \ref{s:main} we show how our main results follow from the Proposition. We emphasize that the concluding argument in Section \ref{s:main} appeals to the -- by now standard -- methods for solving differential inclusions, either by appealing to the Baire category theorem \cite{BressanFlores,Cellina,DacorognaMarcellini97,Bernd01}, or by the more explicit convex integration method \cite{Gromov,MS99,MullerSychev}. 
In our opinion, the Baire category argument
developed in \cite{Bernd} and used in Section \ref{s:main} is,
for the purposes of this paper, the most efficient and
elegant tool.
However, we include in Section \ref{s:CI} an alternative proof 
which follows the convex integration approach, as it makes
easier to "visualize" the solutions constructed in this paper. 

In fact we believe that for $n\geq 3$ a suitable modification of the original approach of Gromov (see \cite{Gromov}) would also work, yielding solutions which are even continuous (work in progress). 

\section{Plane wave analysis of Euler's equations}\label{Tartar}

We start by briefly explaining Tartar's framework \cite{Tartar79}. One considers 
nonlinear PDEs that can be expressed as a system of linear PDEs (conservation laws)
\begin{equation}\label{e:linearpde}
\sum_{i=1}^mA_i\partial_iz=0
\end{equation}
coupled with a pointwise nonlinear constraint (constitutive relations)
\begin{equation}\label{e:pointwiseconst}
z(x)\in K\subset\R^d\textrm{ a.e.,}
\end{equation}
where $z:\Omega\subset\R^m\to\R^d$ is the unknown state variable. The idea is then to consider 
{\it plane wave} solutions to \eqref{e:linearpde}, that is, solutions of the form
\begin{equation}\label{planewave}
z(x)=ah(x\cdot\xi),
\end{equation}
where $h:\R\to\R$. The 
{\it wave cone} $\Lambda$ is given by the states $a\in\R^d$ such that for any choice of the profile $h$
the function \eqref{planewave} solves \eqref{e:linearpde}, that is,
\begin{equation}\label{e:wavecone}
\Lambda:=\left\{a\in\R^d:\,\exists\xi\in\R^m\setminus\{0\}\quad
\mbox{with}\quad \sum_{i=1}^m\xi_i A_ia=0\right\}.
\end{equation}
The oscillatory behavior of solutions to the nonlinear problem is 
then determined by the compatibility of the set $K$ with the 
cone $\Lambda$. 

The Euler equations can be naturally rewritten in this framework. The domain is $\R^m=\R^{n+1}$, and the state variable $z$ is defined as $z=(v,u,q)$, where
$$
q=p+\frac{1}{n}|v|^2,\textrm{ and }u=v\otimes v-\frac{1}{n}|v|^2I_n,
$$
so that $u$ is a symmetric $n\times n$ matrix with vanishing trace and $I_n$ denotes the $n\times n$ identity matrix. From now on the linear space of symmetric $n\times n$ matrices 
will be denoted by $\S^n$ and the subspace of trace--free 
symmetric matrices by $\S^n_0$.
The following lemma is straightforward.

\begin{lemma}\label{LRNC}
Suppose $v\in L^{\infty}(\R^n_x\times\R_t;\R^n)$, $u\in  L^{\infty}(\R^n_x\times\R_t;\S^n_0)$, and 
$q\in L^{\infty}(\R^n_x\times\R_t)$ solve
\begin{equation}\label{LR}
\begin{split}
\partial_tv+\textrm{div }u+\nabla q&=0,\\
\textrm{div }v&=0,
\end{split}
\end{equation}
in the sense of distributions. If in addition
\begin{equation}\label{NC}
u=v\otimes v-\frac{1}{n}|v|^2I_n\quad\textrm{ a.e.~in }\R^n_x\times\R_t,
\end{equation}
then $v$ and $p:=q-\frac{1}{n}|v|^2$ 
are a solution to \eqref{euler} with $f\equiv 0$. Conversely, if $v$ and $p$ solve \eqref{euler} distributionally, then $v$, $u:=v\otimes v-\frac{1}{n}|v|^2I_n$ and 
$q:=p+\frac{1}{n}|v|^2$ solve \eqref{LR} and \eqref{NC}. 
\end{lemma}

Consider the $(n+1)\times(n+1)$ symmetric matrix in block form
\begin{equation}\label{e:U}
U=\begin{pmatrix} u+qI_n&v\\ v&0\end{pmatrix},
\end{equation}
where $I_n$ is the $n\times n$ identity matrix. Notice that by introducing new coordinates $y=(x,t)\in\R^{n+1}$
the equation \eqref{LR} becomes simply
$$
\textrm{div}_yU=0.
$$
Here, as usual, 
a divergence--free matrix field is a matrix of functions with
rows that are divergence--free vectors.
Therefore the wave cone corresponding to \eqref{LR} is given by
$$
\Lambda=\left\{(v,u,q)\in\R^n\times\S^n_0\times\R:\,\det\begin{pmatrix} u+qI_n&v\\ v&0\end{pmatrix}=0\right\}.
$$
\begin{remark}
A simple linear algebra computation shows that for every $v\in\R^n$ and $u\in\S^n_0$ 
there exists $q\in\R$ such that $(v,u,q)\in\Lambda$,
revealing that the wave cone is very large. Indeed, 
let $V^\perp\subset \rn{n}$ be the linear space orthogonal 
to $v$ and consider
on $V^\perp$ the quadratic form $\xi\mapsto \xi\cdot u\xi$.
Then, $\det U = 0$ if and only if $-q$ is an eigenvalue
of this quadratic form.
\end{remark}

In order to exploit this fact for constructing irregular solutions to the nonlinear system, one needs plane wave--like solutions to \eqref{LR} which are localized in space. Clearly an exact plane--wave as in \eqref{planewave} has compact support only if it is identically zero. Therefore this can only be done by introducing an error in the range of the wave, deviating from the line spanned by the wave state $a\in\R^d$. However, this error can be made arbitrarily small. This is the content of the following proposition, which is the building block of our construction.

\begin{proposition}[Localized plane waves]\label{p:lam}
Let $a=(v_0,u_0,q_0)\in \Lambda$ with $v_0\neq 0$, and denote by $\sigma$
the line segment in $\R^n\times\S^n_0\times\R$ joining the points 
$-a$ and $a$.
For every $\eps>0$ there exists a smooth solution
$(v, u, q)$ of \eqref{LR} with the properties:
\begin{itemize}
\item the support of $(v, u, q)$ is contained in $B_1(0)\subset\R^n_x\times\R_t$,
\item the image of $(v,u,q)$ is contained in the $\eps$--neighborhood
of $\sigma$,
\item $\int |v (x,t)|\,dx\,dt \geq \alpha|v_0|$,
\end{itemize}
where $\alpha>0$ is a dimensional constant.
\end{proposition}

\section{Localized plane waves}\label{s:planewaves}

For the proof of Proposition \ref{p:lam} there are two main points. Firstly, we appeal to 
a particular large group of symmetries of the equations in order to reduce the problem to some 
special $\Lambda$-directions. Secondly, to achieve a cut-off which preserves the linear equations \eqref{LR}, we introduce a suitable potential.

\begin{definition}
We denote by $\mathcal{M}$ the set of symmetric $(n+1)\times (n+1)$ matrices $A$
such that $A_{(n+1)(n+1)} = 0$. Clearly, the map
\begin{equation}\label{e:iso}
\R^n\times\S^n_0\times\R\ni (v, u, q) \quad\mapsto\quad
U=\begin{pmatrix} u+qI_n&v\\ v&0\end{pmatrix} \in \mathcal{M}
\end{equation}
is a linear isomorphism. 
\end{definition}
As already observed, in the variables $y=(x, t)\in\R^{n+1}$, the equation \eqref{LR} is equivalent to $\textrm{div }U = 0$.
Therefore Proposition \ref{p:lam} follows immediately 
from 

\begin{proposition}\label{p:lam2}
Let $\overline{U}\in\mathcal{M}$ be
such that $\det\overline{U}=0$ and 
$\overline{U}e_{n+1}\neq 0$, and consider the line segment $\sigma$ with endpoints $-\overline{U}$
and $\overline{U}$. Then there exists 
a constant $\alpha>0$ such that for 
any $\e>0$ there exists a smooth 
divergence--free matrix field 
$U:\R^{n+1}\to\mathcal{M}$ with the properties
\begin{itemize}
\item[(p1)] $\supp U\subset B_1(0)$,
\item[(p2)] $\dist\left(U(y),\sigma\right)<\e$ for all $y\in B_1(0)$,
\item[(p3)] $\int |U(y)e_{n+1}|dy\geq \alpha|\overline{U}e_{n+1}|$,
\end{itemize}
where $\alpha>0$ is a dimensional constant.
\end{proposition}

The proof of Proposition \ref{p:lam2} relies on
two lemmas. The first deals with the 
symmetries of the equations.

\begin{lemma}[The Galilean group]\label{l:transf}
Let $\mathcal{G}$ be the subgroup of $GL_{n+1}(\R)$ defined by
\begin{equation}\label{e:type1}
\big\{A\in \rn{(n+1)\times (n+1)}: \det A \neq 0, Ae_{n+1} = e_{n+1} \big\}\, .
\end{equation}
For every divergence--free map $U: \rn{n+1}\to
\mathcal{M}$ and every $A\in \mathcal{G}$ 
the map 
$$
V (y) := A^t\cdot U (A^{-t} y) \cdot A
$$
is also a divergence--free map $V:\rn{n+1}\to\mathcal{M}$.
\end{lemma}

The second deals with the potential.

\begin{lemma}[Potential in the general case]\label{l:potentialn}
Let $E^{kl}_{ij}\in C^\infty (\R^{n+1})$ be functions for $i,j,k,l =1,\dots, n+1$
so that the tensor $E$ is skew--symmetric in $ij$ and $kl$, that is 
\begin{equation}\label{e:skew}
E^{kl}_{ij} = - E^{lk}_{ij} = -E^{kl}_{ji} = E^{lk}_{ji}\, .
\end{equation}
Then
\begin{equation}\label{e:defU}
U_{ij} = \mathcal{L} (E) = \frac{1}{2} 
\sum_{k,l} \partial^2_{kl} (E^{il}_{kj} + E^{jl}_{ki})
\end{equation}
is symmetric and divergence--free. If in addition 
\begin{equation}\label{e:n+1}
E^{(n+1)j}_{(n+1)i} = 0 \qquad \mbox{for every $i$ and $j$,}
\end{equation}
then $U$ takes values in $\mathcal{M}$.
\end{lemma}

\begin{remark}\label{r:pot2} A suitable potential in the case $n=2$ can be obtained in
a more direct way. Indeed, let $w\in C^{\infty}(\rn{3}, \rn{3})$ be a 
divergence--free vector field and consider the map
$U:\rn{3}\to \mathcal{M}$ given by
\begin{equation}\label{e:potential}
U\;=\; \left(\begin{array}{lll}
\partial_{2} w_1 & \frac{1}{2}\partial_2 w_2 - \frac{1}{2}\partial_{1}
w_1& \frac{1}{2}\partial_{2} w_3\\
\frac{1}{2}\partial_{2} w_2 - \frac{1}{2}\partial_1 w_1&
-\partial_{1} w_2 & - \frac{1}{2}\partial_{1} w_3\\
\frac{1}{2}\partial_{2} w_3& - \frac{1}{2}\partial_{1} w_3& 0
\end{array}\right)\, .
\end{equation}
Then it can be readily checked that 
$U$ is divergence--free. Moreover, $w$ is the curl of
a vector field $\omega$.
However, this is just a particular case of Lemma \ref{l:potentialn}.
Indeed, given $E$ as in the Lemma define the 
tensor $D^k_{ij} = \sum_l \partial_l E^{kl}_{ij}$. Note that $D$ is skew--symmetric in $ij$ and for each $ij$, the vector 
$(D^k_{ij})_{k=1, \ldots, n+1}$ is divergence--free. Moreover,
$$
U_{ij} = \frac{1}{2} \sum_k \partial_k (D_{kj}^i+D_{ki}^j)\, .
$$
Then the vector field $w$ above is simply the special 
choice where $D_{12}^k=-D_{21}^k=w_k$ and all other $D$'s are zero, and a corresponding
relation can be found for $E$ and $\omega$.
\end{remark}

The proofs of the two Lemmas will be postponed until the end of the section and we now
come to the proof of the Proposition.

\begin{proof}[Proof of Proposition \ref{p:lam2}]

\medskip

\noindent{\bf Step 1.} First we treat the case when $\overline{U}\in\mathcal{M}$ is such that 
\begin{equation}\label{e:canonical}
\overline{U}e_1=0,\quad\overline{U}e_{n+1}\neq 0.
\end{equation}
Let
\begin{equation}\label{e:choiceofE}
E_{i1}^{j1}=-E_{1i}^{j1}=-E_{i1}^{1j}=E_{1i}^{1j}=\overline{U}_{ij}\frac{\sin(Ny_1)}{N^2}
\end{equation}
and all the other entries equal to $0$. Note that by our assumption 
$\overline{U}_{ij}=0$ whenever one index is $1$ or both of them
are $n+1$. This ensures that the
tensor $E$ is well defined and satisfies the properties of Lemma \ref{l:potentialn}.

We remark that in the case $n=2$ the matrix $\overline{U}$ takes necessarily the form
\begin{equation}\label{e:canonical2}
\overline{U} \;=\; \left(\begin{array}{lll}
0&0&0\\
0&a&b\\
0&b&0\end{array}\right)\, 
\end{equation} 
with $b\neq 0$,
and we can use the potential of Remark \ref{r:pot2} 
by simply setting 
\begin{equation*}
\begin{split}
w \;&=\; \frac{1}{N} (0, a \cos (N y_1), 2b\cos (N y_1))\,,\\
\omega\; &=\; \frac{1}{N^2} (0, 2b \sin (Ny_1), -a\sin (Ny_1))\, .
\end{split}
\end{equation*}

We come back to the general case. Let $E$ be defined
as in \eqref{e:choiceofE}, fix a smooth cutoff function
$\varphi$ such that
\begin{itemize}
\item $|\varphi|\leq 1$, 
\item $\varphi = 1$ on $B_{1/2}(0)$,
\item $\supp (\varphi)\subset B_1 (0)$, 
\end{itemize}
and consider the map
$$
U = \mathcal{L} (\varphi E).
$$ 
Clearly, $U$ 
is smooth and supported in $B_1 (0)$.
By Lemma \ref{l:potentialn}, 
$U$ is $\mathcal{M}$--valued and 
divergence--free. Moreover 
$$
U(y)=\overline{U}\sin(Ny_1)\textrm{ for }y\in B_{1/2}(0),
$$
and in particular
$$
\int |U(y)e_{n+1}|dy\geq |\overline{U}e_{n+1}| 
\int_{B_{1/2}(0)}|\sin (Ny_1)|\, dy \geq 2\alpha|\overline{U}e_{n+1}|,
$$ 
for some positive dimensional constant $\alpha=\alpha(n)$ for sufficiently large $N$. 

Finally, observe that 
$$
U - \varphi\tilde{U} \;=\; \mathcal{L} 
(\varphi E)
- \varphi \mathcal{L} (E)
$$
is a sum of products of first--order 
derivatives of $\varphi$ with
first--order derivatives of components 
of $E$ and of second--order derivatives of $\varphi$ 
with components of $E$. 
Thus, 
$$
\|U - \varphi \tilde{U}\|_\infty \;\leq\;
C \|\varphi\|_{C^2} \|E\|_{C^1} \;\leq\; 
\frac{C'}{N} \|\varphi\|_{C^2}\, ,
$$ 
and by choosing $N$ sufficiently large we obtain 
$\|U-\varphi \tilde{U}\|_{\infty}< \eps$. On the other hand, 
since $|\varphi|\leq 1$ and $\tilde{U}$ takes values in $\sigma$, 
the image of $\varphi\tilde{U}$ is also contained in 
$\sigma$. This shows that the image of $U$ is contained in
the $\eps$--neighborhood of $\sigma$.

\medskip

\noindent{\bf Step 2.} We treat the general case by reducing to the situation above. Let $\overline{U}\in\mathcal{M}$ be as in the Proposition, so that
$$
\overline{U}f=0,\quad \overline{U}e_{n+1}\neq 0,
$$
where $f\in\R^{n+1}\setminus\{0\}$ is such that $\{f,e_{n+1}\}$ are linearly independent. 
Let $f_1, \ldots, f_{n+1}$ be a basis for $\R^{n+1}$ such that $f_1 = f$ and $f_{n+1}
= e_{n+1}$ and consider the matrix $A$ such that 
$$
A e_i = f_i\textrm{ for }i=1,\dots,n+1.
$$ 
Then $A\in \mathcal{G}$ 
(cf.~with the definition of $\mathcal{G}$ given 
in Lemma \ref{l:transf}), 
and the map
\begin{equation}\label{e:linearisom}
T:X\mapsto (A^{-1})^tXA^{-1}
\end{equation}
is a linear isomorphism of $\R^{n+1}$.
Set 
\begin{equation}\label{e:specialcase}
\overline{V} = A^t\overline{U} A,
\end{equation}
so that $\overline{V}\in \mathcal{M}$ satisfies
$$
\overline{V}e_1=0,\quad\overline{V}e_{n+1}\neq 0.
$$
Given $\eps>0$, using Step 1 we construct a smooth map $V:\R^{n+1}\to\mathcal{M}$ 
supported in $B_1 (0)$
with the image lying in the $\|T\|^{-1}\eps$--neighborhood 
of the line segment $\tau$ with endpoints
$-\overline{V}$ and $\overline{V}$, and such that
$$
V(y)=\overline{V}\sin(Ny_1).
$$ 
Let $U$ be the 
$\mathcal{M}$--valued map 
$$
U (y) = (A^{-1})^t V(A^ty) A^{-1}.
$$ 
By our discussion above the isomorphism
$T:X\mapsto (A^{-1})^tXA^{-1}$ maps the line segment $\tau$ onto $\sigma$.
Therefore:
\begin{itemize}
\item $U$ is supported in $A^{-t} (B_1 (0))$ and it is smooth,
\item $U$ is divergence--free thanks to Lemma \ref{l:transf},
\item $U$ takes values in an $\eps$--neighborhood of
the segment $\sigma$,
\end{itemize}
and furthermore
\begin{eqnarray}
\int_{A^{-t} (B_1 (0))} |U(y)e_{n+1}| dy
&=& \int_{A^{-t} (B_1(0))}|A^{-t}V(A^ty)e_{n+1}|dy\nonumber\\
&=&\int_{B_1 (0)} |A^{-t} V (z) e_{n+1}| \frac{dz}{|\det A^t|}
\nonumber\\
&\geq& \frac{2\alpha|A^{-t}\overline{V}e_{n+1}|}{|\det A|}
\;=\; \frac{2\alpha}{|\det A|} |\overline{U}e_{n+1}|.
\label{e:changevar}
\end{eqnarray}
To complete the proof we appeal to a standard covering/rescaling argument. That is, we can find a finite number of points $y_k\in B_1(0)$ and radii $r_k>0$ so that the rescaled and 
translated sets $A^{-t} (B_{r_k}(y_k))$ 
are pairwise disjoint, all contained in $B_1(0)$, and
\begin{equation}\label{e:covers}
\bigcup_{k} \left|A^{-t}(B_{r_k}(y_k))\right|
\geq \frac{1}{2}|B_1(0)|.
\end{equation}
Let $U_k(y)=U(\frac{y-y_k}{r_k})$ and $\tilde U=\sum_k U_k$. Then $\tilde U:\R^{n+1}\to\mathcal{M}$ is smooth, clearly satisfies (p1) and (p2), and
\begin{eqnarray}
\int|\tilde U(y)e_{n+1}|dy&=&
\sum_k\int_{A^{-t}B_{r_k}(y_k)}|U_k(y)e_{n+1}|dy\nonumber\\
&\stackrel{\eqref{e:changevar}}{\geq}&
\sum_k 2\alpha|\overline{U}e_{n+1}| |\det A|^{-1} 
\frac{|B_{r_k}(y_k)|}{|B_1(0)|}\nonumber\\
&=& 2\alpha |\overline{U}e_{n+1}|
\frac{\sum_k \left|A^{-t} (B_{r_k} (y_k))\right|}{|B_1 (0)|}
\;\stackrel{\eqref{e:covers}}{\geq}\; \alpha|\overline{U}e_{n+1}|.
\nonumber\end{eqnarray}
This completes the proof.
\end{proof} 

\begin{proof}[Proof of Lemma \ref{l:transf}]
First of all we check that whenever $B\in\mathcal{M}$, then $A^tBA\in\mathcal{M}$ for 
all $A\in\mathcal{G}$.
Indeed, $A^tBA$ is symmetric, and 
since $A$ satisfies $Ae_{n+1}=e_{n+1}$, we have
\begin{eqnarray}
(A^t BA)_{(n+1)(n+1)} &=& e_{n+1}\cdot A^tBAe_{n+1}=Ae_{n+1}\cdot BAe_{n+1}\nonumber\\
&=&e_{n+1}\cdot Be_{n+1}=B_{(n+1)(n+1)}=0. \label{e:straight1}
\end{eqnarray}

Now, let $A$, $U$ and $V$ be as in the statement. The argument
above shows that $V$ is $\mathcal{M}$--valued. It remains to check
that if $U$ is divergence--free, then $V$ is also divergence--free.
To this end let $\phi\in C^{\infty}_c(\R^{n+1};\R^{n+1})$ be a compactly supported test function
and consider $\tilde\phi\in C^{\infty}_c(\R^{n+1};\R^{n+1})$ defined by 
$$
\tilde\phi(x)=A\phi(A^tx).
$$
Then $\nabla\tilde\phi(x)=A\nabla\phi(A^tx)A^t$, and
by a change of variables we obtain
\begin{equation*}
\begin{split}
\int \rm tr\bigl(V(y)\nabla\phi(y)\bigr)dy=&\int\rm tr\bigl(A^tU(A^{-t}y)A\nabla\phi(y)\bigr)dy\\
=&\int\rm tr\bigl(U(A^{-t}y)A\nabla\phi(y)A^t\bigr)dy\\
=&\int\rm tr\bigl(U(x)A\nabla\phi(A^tx)A^t\bigr)(\det A)^{-1}dx\\
=&(\det A)^{-1}\int\rm tr\bigl(U(x)\nabla\tilde\phi(x)\bigr)dx=0,
\end{split}
\end{equation*}
since $U$ is divergence--free. But this implies that $V$ is also divergence-free.

\end{proof}

\begin{proof}[Proof of Lemma \ref{l:potentialn}]
First of all, $U$ is clearly symmetric and $U_{(n+1)(n+1)}=0$. 
Hence $U$ takes values in $\mathcal{M}$. 
To see that $U$ is divergence--free, we calculate
\begin{eqnarray*}
\sum_j\partial_jU_{ij}&=&
 \frac{1}{2} 
\sum_{k,l} \partial^3_{jkl} (E^{il}_{kj} + E^{jl}_{ki})\\
&=& \frac{1}{2}\sum_l \partial_l \Big(\sum_{jk} \partial^2_{jk}
E^{il}_{kj}\Big) + \frac{1}{2} \sum_k \partial_k \Big( \sum_{jl}
\partial^2_{jl} E^{jl}_{ki}\Big)\;\stackrel{\eqref{e:skew}}{=}\; 0\, .
\end{eqnarray*}
This completes the proof of the lemma.
\end{proof}

%%%%%%%%%%%%%%%%%%%%%%%%%%%%%%%%%%%

\section{Proof of the main results}\label{s:main}

For clarity we now state the precise form of our main result. Theorems \ref{main},  \ref{t:decrease}  and \ref{main2} are direct corollaries.

\begin{theorem}\label{t:precise}
Let $\Omega\subset\R^n_x\times\R_t$ be a bounded open domain. There exists $(v,p)\in L^{\infty}(\R^n_x\times\R_t)$ solving the Euler equations
\begin{equation*}
\begin{split}
\partial_tv+\div(v\otimes v)+\nabla p &= 0\\
\div v&=0\, ,
\end{split}
\end{equation*}
such that 
\begin{itemize}
\item $|v(x,t)|=1$ for a.e.~$(x,t)\in\Omega$,
\item $v(x,t)=0$ and $p(x,t)=0$ for a.e.~$(x,t)\in(\R^n_x\times\R_t)\setminus\Omega$.
\end{itemize}
Moreover, there exists a sequence of functions $(v_k,p_k,f_k)\in C^{\infty}_c(\Omega)$ such that
\begin{equation*}
\begin{split}
\partial_tv_k+\div(v_k\otimes v_k)+\nabla p_k &= f_k\\
\div v_k&=0\, ,
\end{split}
\end{equation*}
and
\begin{itemize}
\item $f_k$ converges to $0$ in $H^{-1}$,
\item $\|v_k\|_\infty+ \|p_k\|_\infty$ 
is uniformly bounded,
\item $(v_k,p_k)\to (v,p)$
in $L^q$ for every $q<\infty$.
\end{itemize}
\end{theorem}

We remark that the statements of Theorem \ref{main} and
Theorem \ref{main2} are just subsets of the statement of
Theorem \ref{t:precise}. As for Theorem \ref{t:decrease},
note that it suffices to choose, for instance, $\Omega =
B_r (0) \times ]-1,1[$, where $B_r (0)$ is the ball of $\R^n$
with volume $1$. 

We recall from Lemma \ref{LRNC} that for the first half of the theorem it suffices to prove that 
there exist
$$
(v,u,q)\in L^{\infty}(\R^n_x\times\R_t;\R^n\times\S^n_0\times\R)
$$ 
with support in $\Omega$, such that $|v|=1$ a.e.~in $\Omega$ and 
(\ref{LR}) and (\ref{NC}) are satisfied. In Proposition 
\ref{p:lam} we constructed compactly supported solutions 
$(v,u,q)$ to (\ref{LR}). The point is thus to find solutions 
which satisfy in addition the pointwise constraint (\ref{NC}). The main idea is to consider
the sets
\begin{equation}\label{e:Kdef}
K=\left\{(v,u)\in \R^n\times\S^n_0:\,u=v\otimes v-\frac{1}{n}|v|^2I_n\,,\; |v|=1\right\},
\end{equation}
and 
\begin{equation}\label{e:Udef}
\mathcal{U}=\textrm{int }(K^{co}\times [-1,1]),
\end{equation}
where $\textrm{int }$ denotes the topological interior of the set in $\R^n\times\S^n_0
\times\R$, and $K^{co}$ denotes the convex hull of $K$.
Thus, a triple $(v,u,q)$ solving \eqref{LR} and taking values in the convex extremal points of $\overline{\mathcal{U}}$ is indeed a solution to \eqref{NC}. We will prove that $0\in\mathcal{U}$, and therefore there exist plane waves taking values in $\mathcal{U}$. 
The goal is to add them so to get an infinite sum 
$$
(v,u,q)=\sum_{i=1}^{\infty}(v_i,u_i,q_i)
$$ 
with the properties that
\begin{itemize}
\item the partial sums $\sum_{i=0}^k (v_i, u_i, q_i)$ take values in $\mathcal{U}$,
\item $(v,u,q)$ is supported in $\Omega$,
\item $(v,u,q)$ takes values in the convex extremal points of $\overline{U}$ a.e.~in $\Omega$,
\item $(v,u,q)$ solves the linear partial differential equations \eqref{LR}.
\end{itemize}
There are two important reasons why this construction is possible.
First of all, since the wave cone $\Lambda$ is very large, we can always get closer and closer to the extremal point of $\mathcal{U}$ with the sequence $(v_k, u_k, p_k)$. Secondly, because the waves are localized in space--time, by choosing the supports smaller and smaller we can achieve 
strong convergence of the sequence. 
In view of Lemma \ref{LRNC} this gives the solution of Euler that we we are looking for. The  
partial sums give the approximating sequence of the theorem. 

\medskip

This sketch of the proof is philosophically closer to the method of convex integration, where the difficulty is to ensure strong convergence of the partial sums. The Baire category argument avoids this difficulty by introducing a metric for the space of solutions to \eqref{LR} with values in $\mathcal{U}$, and proving that in its closure a generic element takes values in the convex extreme points. An interesting corollary of the Baire category argument is that, within the class of solutions to the Euler equations with driving force in some particular bounded subset of $H^{-1}$, the typical (in the sense of category) element has the properties of Theorem \ref{t:precise} .

We split the proof of Theorem \ref{t:precise} 
into several lemmas and 
a short concluding argument, which is given at the
beginning of Section \ref{s:cont}. For the purpose
of this section, we could have
presented a shorter proof, avoiding Lemma
\ref{l:geometric} and without giving the explicit bound
\eqref{e:gain} of Lemma \ref{l:onestep}. However,
these statements will be needed in the convex integration
proof of Section \ref{s:CI}.

\subsection{The geometric setup}

\begin{lemma}\label{l:convexhull}
Let $K$ and $\mathcal{U}$ be defined as in \eqref{e:Kdef} and \eqref{e:Udef}, i.e.~
\begin{equation*}
K=\left\{(v,u)\in\mathbb{S}^{n-1}\times\S^n_0:\,u=v\otimes v-\frac{I_n}{n}\right\}.
\end{equation*}
Then $0\in\textrm{int }K^{co}$ and hence $0\in\mathcal{U}$.
\end{lemma}

\begin{proof}
Let $\mu$ be the Haar measure on $\mathbb{S}^{n-1}$ and consider the linear map
$$
T:C(\mathbb{S}^{n-1})\to\R^n\times \S^n_0,\quad \phi\mapsto 
\int_{\mathbb{S}^{n-1}}\left(v,v\otimes v-\frac{I_n}{n}\right)
\phi(v)\,d\mu\, .
$$
Clearly, if 
\begin{equation}\label{e:convexity}
\phi\geq 0  \qquad\mbox{and} \qquad
\int_{\mathbb{S}^{n-1}} \phi\, d\mu =1\, ,
\end{equation} 
then $T (\phi)\in K^{co}$.
Notice that
$$
T (1) \;=\; \int_{\mathbb{S}^{n-1}}\left(v,v\otimes v-\frac{I_n}{n}\right)\,d\mu\;=\;0,
$$
and hence $0\in K^{co}$. Moreover, whenever $\psi\in C (\mathbb{S}^{n-1})$ is such that
\begin{equation}\label{e:convexity2}
\alpha \;=\; 1 - \int_{\mathbb{S}^{n-1}} \psi\, d\mu\;\geq\; 
\|\psi\|_{C ({\mathbb{S}^{n-1}})},
\end{equation}
$\phi = \alpha +\psi$ satisfies \eqref{e:convexity} and hence $T (\psi) = T (\phi)\in
K^{co}$. Since \eqref{e:convexity2} holds whenever $\|\psi\|_{C ({\mathbb{S}^{n-1}})}<1/2$,
it suffices to show that $T$ is surjective to prove that $K^{co}$ contains a neighborhood
of $0$.

The surjectivity of $T$ follows from orthogonality in $L^2(\mathbb{S}^{n-1})$. 
Indeed, letting $\phi=v_i$ for each $i$, we obtain 
$$
T(\phi)=\beta_1(e_i,0)\textrm{, where }\beta_1=\int_{\mathbb{S}^{n-1}}v_1^2d\mu.
$$ 
Furthermore, setting $\phi=v_iv_j$ with $i\neq j$, we obtain 
$$
T(\phi)=\beta_2\bigl(0,e_i\otimes e_j+e_j\otimes e_i\bigr)\textrm{, where }\beta_2=\int_{\mathbb{S}^{n-1}}v_1^2v_2^2d\mu.
$$ 
Finally, setting 
$\phi=v_i^2-\frac{1}{n}$ we obtain
$$
T(\phi)=\beta_3\Bigl(0,e_i\otimes e_i-\frac{1}{(n-1)}\sum_{j\neq i}e_j\otimes e_j\Bigr),
$$
where 
$$
\beta_3=\int_{\mathbb{S}^{n-1}}\left(v_1^2-\frac{1}{n}\right)^2d\mu.
$$ 
This shows that the image of $T$ contains $n+\frac{1}{2}n(n+1)-1$ linearly independent elements, hence a basis for $\R^n\times\S^n_0$.

\end{proof}

\bigskip

\begin{lemma}\label{l:geometric}
There exists a dimensional constant $C>0$ such that for any 
$(v,u,q)\in \mathcal{U}$ there exists 
$(\bar{v},\bar{u})\in\R^n\times\S^n_0$ such that
$(\bar{v},\bar{u},0)\in\Lambda$,
the line segment with endpoints $(v,u,q)\pm(\bar{v},\bar{u},0)$ is contained in 
$\mathcal{U}$, and 
$$
|\bar{v}|\geq C(1-|v|^2).
$$
\end{lemma}

\begin{proof}
Let $z=(v,u)\in \textrm{int }K^{co}$. By Carath\'eodory's theorem $(v,u)$ lies in the interior of a simplex in $\R^n\times\S^n_0$ spanned by elements of $K$. In other words
$$
z=\sum_{i=1}^{N+1}\lambda_iz_i,
$$
where $\lambda_i\in\, ]0,1[\,$, $z_i=(v_i,u_i)\in K$, 
$\sum_{i=1}^{N+1}\lambda_i=1$, 
and $N=n(n+3)/2-1$ is the dimension of $\R^n\times\S^n_0$. Assume that the coefficients are ordered so that $\lambda_1=\max_i\lambda_i$. Then for any $j>1$
$$
z\,\pm\,\frac{1}{2}\lambda_j(z_j-z_1)\in \textrm{int }K^{co}.
$$
Indeed,
\begin{equation*}
z\,\pm\,\frac{1}{2}\lambda_j(z_j-z_1)=\sum_i\mu_iz_i,
\end{equation*}
where $\mu_1=\lambda_1\mp \frac{1}{2}\lambda_j$, $\mu_j=\lambda_j\pm\frac{1}{2}\lambda_j$ and $\mu_i=\lambda_i$ for $i\notin\{1,j\}$. It is easy to see that $\mu_i\in\, ]0,1[$ for all $i=1\dots N+1$. 

On the other hand $z-z_1=\sum_{i=2}^{N+1}\lambda_i(z_i-z_1)$, so that in particular
$$
|v-v_1|\leq N\max_{i=2\dots N+1}\lambda_i|v_i-v_1|
$$
Let $j>1$ be such that $\lambda_j|v_j-v_1|=\max_{i=2\dots N+1}\lambda_i|v_i-v_1|$, and
let 
$$
(\bar{v},\bar{u})=\frac{1}{2}\lambda_j(z_j-z_1).
$$
The line segment with endpoints 
$(v,u)\pm(\bar{v},\bar{u})$ is contained in the interior of $K^{co}$ and
hence also the line segment $(v,u,q)\pm(\bar{v},\bar{u},0)$ is contained in $\mathcal{U}$. Furthermore
$$
\frac{1}{4N}(1-|v|^2)\leq \frac{1}{2N}(1-|v|)\leq \frac{1}{2N}(|v-v_1|)\leq |\bar{v}|.
$$
Finally, we show that $(\bar{v},\bar{u},0)\in\Lambda$. This amounts to showing
that whenever $a,b\in \mathbb{S}^{n-1}$, the matrix
$$
\begin{pmatrix}a\otimes a-\frac{I_n}{n}&a\\ a&0\end{pmatrix}-
\begin{pmatrix}b\otimes b-\frac{I_n}{n}&b\\ b&0\end{pmatrix}
$$
has zero determinant and hence lies in the wave cone $\Lambda$ defined in 
\eqref{e:wavecone}. Let $P\in GL_n(\R)$ with $Pa=e_1$ and $Pb=e_2$. Note that
$$
\begin{pmatrix}P&0\\ 0&1\end{pmatrix} \begin{pmatrix}a\otimes a&a\\ a&0\end{pmatrix}
\begin{pmatrix}P^t&0\\ 0&1\end{pmatrix}=\begin{pmatrix}Pa\otimes Pa&Pa\\ Pa&1\end{pmatrix},
$$
so that it suffices to check the determinant of 
$$
\begin{pmatrix}e_1\otimes e_1&e_1\\ e_1&0\end{pmatrix}-
\begin{pmatrix}e_2\otimes e_2&e_2\\ e_2&0\end{pmatrix}.
$$
Since $e_1+e_2-e_{n+1}$ is in the kernel of this matrix, it has indeed determinant zero. This completes the proof.

\end{proof}

\subsection{The functional setup}

We define the complete metric space $X$ as follows. Let
$$
X_0:=\bigl\{(v,u,q)\in C^{\infty}(\R^n_x\times\R_t):\;
\mbox{(i), (ii) and (iii) below hold}       
\bigr\} 
$$
\begin{itemize}
\item[(i)] $\supp(v,u,q)\subset \Omega$,
\item[(ii)] $(v,u,q)$ solves \eqref{LR} in $\R_x^n\times\R_t$,
\item[(iii)] $(v (x,t),u(x,t),q(x,t)) \in\mathcal{U}$ 
for all $(x,t)\in\R^n_x\times\R_t$.
\end{itemize}
We equip $X_0$ with the topology of $L^{\infty}$-weak* 
convergence of $(v,u,q)$ and we let
$X$ be the closure of $X_0$ in this topology. 

\begin{lemma}\label{l:metricspace}
The set $X$ with the topology of $L^{\infty}$ weak* convergence is a nonempty compact metrizable space. Moreover, if $(v,u,q)\in X$ is such that
$$
|v(x,t)|=1\quad\textrm{ for almost every }(x,t)\in\Omega,
$$
then $v$ and $p:=q-\frac{1}{n}|v|^2$ is a weak solution of \eqref{euler} in $\R^n_x\times\R_t$ such that
$$
v(x,t)=0\textrm{ and }p(x,t)=0\textrm{ for all }(x,t)\in\R^n_x\times\R_t\setminus\Omega.
$$
\end{lemma}

\begin{proof}
In Lemma \ref{l:convexhull} we showed that $0\in\mathcal{U}$, hence $X$ is nonempty. Moreover, $X$ is a bounded and closed subset of $L^{\infty}(\Omega)$,
hence with the weak* topology it becomes a compact metrizable space. Since $\overline{\mathcal{U}}$ is a compact convex set, any $(v,u,q)\in X$ satisfies
$$
\supp(v,u,q)\subset\overline{\Omega},\quad
(v,u,q) \textrm{ solves \eqref{LR} and takes values in }\overline{\mathcal{U}}.
$$
In particular $(v,u)(x,t)\in K^{co}$ almost everywhere. 
Finally, observe also that if $(v,u)(x,t)\in K^{co}$, then
$$
(v,u)(x,t)\in K\textrm{ if and only if }|v(x,t)|=1.
$$
In light of Lemma \ref{LRNC} this concludes the proof.
\end{proof}

Fix a metric $d^*_{\infty}$ inducing the weak* topology of $L^{\infty}$ in $X$, so that $(X,d^*_{\infty})$ is a complete metric space. 

\begin{lemma}\label{l:identitymap}
The identity map
$$
I: (X,d^*_{\infty})\to L^2(\R^n_x\times\R_t)\,\textrm{ defined by }(v,u,q)\mapsto (v,u,q)
$$
is a Baire-1 map and therefore the set of points of continuity is residual in $(X,d^*_{\infty})$.
\end{lemma}
\begin{proof}
Let $\phi_r(x,t)=r^{-(n+1)}\phi(rx,rt)$ be any regular space-time convolution kernel. 
For each fixed $(v,u,q)\in X$ we have
$$
(\phi_r*v,\phi_r*u,\phi_r*q)\to (v,u,q)\textrm{ strongly in }L^2\textrm{ as }r\to 0.
$$
On the other hand, for each $r>0$ and $(v^k,u^k,q^k)\in X$ 
$$
(v^k,u^k,q^k)\overset{*}{\rightharpoonup}(v,u,q)\textrm{ in }L^{\infty}\,\Longrightarrow\, 
\phi_r*(v^k,u^k,q^k)\to\phi_r*(v,u,q)\textrm{ in }L^2.
$$
Therefore each map $I_r: (X,d^*_{\infty})\to L^2$ defined by 
$$
I_r:(u,v,q)\mapsto (\phi_r*v,\phi_r*u,\phi_r*q)
$$
is continuous, and 
$$
I(v,u,q)=\lim_{r\to 0}I_r(v,u,q)\qquad
\textrm{ for all } (v,u,q)\in X\, .
$$ 
This shows that $I:(X,d^*_{\infty})\to L^2$ is a pointwise limit of continuous maps, hence it is a Baire-1 map. Therefore the set of points of continuity of $I$ is residual in $(X,d^*_{\infty})$, see \cite{Oxtoby}. 
\end{proof}

\bigskip

\subsection{Points of continuity of the identity map}\label{s:cont}

The proof of Theorem \ref{t:precise} will follow from Lemmas \ref{l:metricspace}
and \ref{l:identitymap} once we prove the following 

\bigskip

\noindent{\bf Claim: }
If $(v,u,q)\in X$ is a point of continuity of $I$, 
then 
\begin{equation}\label{e:claim}
|v(x,t)|=1 \textrm{ for almost every }(x,t)\in\Omega\, .
\end{equation}

\bigskip

Indeed, if the claim is true, then the set of $(v,u,q)\in X$ such that
$|v|=1$ a.e.~is nonempty, yielding solutions of \eqref{euler}. Furthermore,
any such $(v,u,q)$ must be the 
strong $L^2$ limit of some sequence $\{(v_k,u_k,q_k)\}
\subset X_0$. Therefore, with $p_k=q_k-\frac{1}{n}|v_k|^2$, and 
$$
f_k=\textrm{div }\left(v_k\otimes v_k-\frac{1}{n}|v_k|^2Id-u_k\right),
$$
we obtain $\textrm{div }v_k=0$ and
$$
\partial_tv_k+\textrm{div }v_k\otimes v_k+\nabla p_k=f_k.
$$
Moreover, $f_k\to 0$ in $H^{-1}$. 

\bigskip

Therefore it remains to prove our claim. Observe that since
$|v(x,t)|\leq 1$ a.e.~$(x,t)\in\Omega$, \eqref{e:claim} is equivalent to
$$
\|v\|_{L^2(\Omega)}=|\Omega|,
$$
where $|\Omega|$ denotes the $(n+1)$-dimensional Lebesgue measure of $\Omega$.
To prove the claim we prove the following lemma, from which the claim immediately follows:

\begin{lemma}\label{l:onestep}
There exists a dimensional constant $\beta>0$ with the following property.
Given $(v_0,u_0,q_0)\in X_0$ there exists a sequence 
$(v_k,u_k,q_k)\in X_0$
such that
\begin{equation}\label{e:gain}
\|v_k\|^2_{L^2(\Omega)}\geq \|v_0\|^2_{L^2(\Omega)}+\beta
\left(|\Omega|-\|v_0\|^2_{L^2(\Omega)}\right)^2,
\end{equation}
and
$$
(v_k,u_k,q_k)\overset{*}{\rightharpoonup} (v_0,u_0,q_0)\textrm{ in }L^{\infty}(\Omega).
$$
\end{lemma}

Indeed, assume for a moment that $(v,u,q)$ is a point
of continuity of $I$. Fix a sequence
$\{(v_k, u_k, q_k\}\subset X_0$ converges weakly$^*$ to
$(v,u,q)$. Using Lemma \ref{l:onestep} and a standard
diagonal argument, we can produce a second sequences
$(\tilde{v}_k, \tilde{u}_k, \tilde{q}_k)$ which converges
weakly$^*$ to $(v,u,q)$ and such that
\begin{equation}\label{l:liminf}
\liminf_{k\to\infty} \|\tilde{v}_k\|_2^2 
\;\geq\; \liminf_{k\to\infty} \left(\|v_k\|^2_2
+ \beta \left(|\Omega| - \|v_k\|_2^2\right)^2\right)\, .
\end{equation}
Since $I$ is continuous at $(v,u,q)$, both $v_k$ and
$\tilde{v}_k$ converge strongly to $v$. Therefore
\begin{equation}
\|v\|_2^2 \;\geq\; \|v\|_2^2 + \beta
\left(|\Omega| - \|v\|_2^2\right)^2\, .
\end{equation}
Therefore, $\|v\|^2_2 = |\Omega|$. On the other hand,
since $v=0$ a.e. outside $\Omega$ and $|v|\leq 1$ a.e.
on $\Omega$, this implies \eqref{e:claim}.
 
\begin{proof}[Proof of Lemma \ref{l:onestep}]

\noindent{\bf Step 1. }
Let $(v_0,u_0,q_0)\in X_0$. By Lemma \ref{l:geometric}
for any $(x,t)\in\Omega$ there exists a direction 
$$
\bigl(\bar{v}(x,t),\bar{u}(x,t)\bigr)\in\R^n\times\S^n_0
$$
such that the line segment with endpoints
$$
\bigl(v_0(x,t),u_0(x,t),q_0(x,t)\bigr)\pm\bigl(\bar{v}(x,t),\bar{u}(x,t),0\bigr)
$$
is contained in $\mathcal{U}$, and
$$
|\bar{v}(x,t)|\geq C(1-|v_0(x,t)|^2).
$$
Moreover, since $(v_0,u_0,q_0)$ is uniformly continuous, there exists $\eps>0$ such that
for any $(x,t),(x_0,t_0)\in \Omega$ with $|x-x_0|+|t-t_0|<\eps$, the $\eps$-neighbourhood of 
the line segment with endpoints
$$
\bigl(v_0(x,t),u_0(x,t),q_0(x,t)\bigr)\pm\bigl(\bar{v}(x_0,t_0),\bar{u}(x_0,t_0),0\bigr)
$$
is also contained in $\mathcal{U}$.

\medskip

\noindent{\bf Step 2. } 
Fix $(x_0,t_0)\in\Omega$ for the moment. Use Proposition \ref{p:lam} with
$$
a=(\bar{v}(x_0,t_0),\bar{u}(x_0,t_0),0)\in\Lambda
$$
and $\eps>0$ to 
obtain a smooth solution $(v,u,q)$ of \eqref{LR} with the properties
stated in the Proposition, and for any $r<\eps$ let
$$
(v_r,u_r,q_r)(x,t)=(v,u,q)\left(
\frac{x-x_0}{r},\frac{t-t_0}{r}\right).
$$
Then $(v_r,u_r,q_r)$ is also a smooth solution of \eqref{LR}, with the properties
\begin{itemize}
\item the support of $(v_r, u_r, q_r)$ is contained in $B_r(x_0,t_0)\subset\R^n_x\times\R_t$,
\item the image of $(v_r,u_r,q_r)$ is contained in the $\eps$--neighborhood
of the line-segment with endpoints $\pm(\bar{v}(x,t),\bar{u}(x,t),0)$,
\item and 
$$
\int |v_r (x,t)|\,dx\,dt \geq \alpha|\bar{v}(x_0,t_0)||B_r(x_0,r_0)|.
$$
\end{itemize}
In particular, for any $r<\eps$ we have
$(v_0,u_0,q_0)+(v_r,u_r,q_r)\in X_0$.

\medskip

\noindent{\bf Step 3. }
Next, observe that since $v_0$ is uniformly continuous, there exists $r_0>0$ such that
for any $r<r_0$ there exists a finite family of pairwise disjoint balls $B_{r_j}(x_j,t_j)\subset\Omega$ with $r_j<r$ such that
\begin{equation}\label{e:CI-discretize}
\int_{\Omega}(1-|v_0(x,t)|^2)dxdt\leq 2\sum_{j}(1-|v_0(x_j,t_j)|^2)|B_r(x_j,t_j)|
\end{equation}

Fix $k\in\N$ with $\frac{1}{k}<\min\{r_0,\eps\}$ and choose a finite family of pairwise disjoint 
balls $B_{r_{k,j}}(x_{k,j},t_{k,j})\subset\Omega$ with radii $r_{k,j}<\frac{1}{k}$ such that \eqref{e:CI-discretize} holds.
In each ball $B_{r_{k,j}}(x_{k,j},t_{k,j})$ we apply the construction above to obtain
$(v_{k,j},u_{k,j},q_{k,j})$, and in particular we then have
$$
(v_k,u_k,q_k):=(v_0,u_0,q_0)+\sum_j(v_{k,j},u_{k,j},q_{k,j})\in X_0,
$$
and
\begin{eqnarray}
\int |v_k(x,t)-v_0(x,t)|dxdt&=&
\sum_j\int |v_{k,j}(x,t)|dxdt\nonumber\\
&\geq& \alpha\sum_j|\bar{v}(x_{k,j},t_{k,j})|
|B_{r_{k,j}}(x_{k,j},t_{k,j})|\nonumber\\
&\geq& C\alpha\sum_j
(1-|v_0(x_{k,j},t_{k,j})|^2)|B_{r_{k,j}}(x_{k,j},t_{k,j})|
\nonumber\\
&\geq& \frac{1}{2}C\alpha\int_{\Omega}(1-|v_0(x,t)|^2)dxdt.
\label{e:L1est}
\end{eqnarray}
Finally observe that by letting $k\to \infty$, the above construction yields a sequence
$(v_k,u_k,q_k)\in X_0$ such that
\begin{equation}\label{e:weaks}
(v_k,u_k,q_k)\overset{*}{\rightharpoonup} (v_0,u_0,q_0).
\end{equation}
Hence,
\begin{eqnarray}
\liminf_{k\to \infty}\|v_k\|_{L^2 (\Omega)}&=&
\|v_0\|^2_2+
\liminf_{k\to\infty} \left(\langle v_0,
(v_k - v_0)\rangle_2
+ \|v_k-v_0\|^2_2\right)\nonumber\\
&\stackrel{\eqref{e:weaks}}{=}&
|v_0\|^2_2 + \liminf_{k\to\infty}
\|v_k-v_0\|^2_2\nonumber\\
&\geq&\|v_0\|^2_2 + |\Omega| \liminf_{k\to\infty}
\left(\|v_k-v_0\|_{L^1 (\Omega)}\right)^2 \label{e:uno}
\end{eqnarray}
Combining \eqref{e:L1est} and \eqref{e:uno} we
get
\begin{eqnarray}
\liminf_{k\to \infty}\|v_k\|_{L^2 (\Omega)}
&\geq& 
\|v_0\|^2_{L^2(\Omega)}+
\frac{|\Omega| C^2\alpha^2}{4}
\left(|\Omega|-\|v_0\|^2_{L^2(\Omega)}\right)^2\nonumber,
\end{eqnarray}
which gives \eqref{e:gain} with
$\beta=\frac{1}{4}|\Omega| C^2\alpha^2$.
\end{proof}

\section{A proof of Theorem \ref{t:precise} using convex
integration}\label{s:CI}

In this section we provide an alternative, more direct proof for Theorem \ref{t:precise}, following the method of convex integration as presented for example in \cite{MS99}. 

In fact the two approaches (i.e. Baire category methods and convex integration) can be unified to a large extent. For a discussion comparing the two approaches we refer to the end of Section 3.3 in \cite{Bernd}, see also the paper \cite{Sychev-FR} for a different point of view. Nevertheless, in order to get a feeling for the type of solution that Theorem \ref{t:precise} produces, it helps to see the direct construction of the convex integration method. 

We will freely refer to the notation of the previous sections. In particular the proof relies on Lemmas \ref{l:convexhull}, \ref{l:geometric}, \ref{l:metricspace} and \ref{l:onestep}. These results enable us to construct an approximating sequence, as explained briefly at the beginning of Section \ref{s:main}, by adding (almost-)plane-waves on top of each other. It is only the limiting step that is more explicit in this approach. 
The following argument is essentially from Section 3.3 in \cite{MS99}.

\begin{proof}[Alternative proof of Theorem \ref{t:precise}]
Using Lemma \ref{l:onestep}, we construct inductively a 
sequence $(v_k,u_k,q_k)\in X_0$ and a sequence of numbers $\eta_k>0$ as follows. 
Let $\rho_{\eps}$ be a standard mollifying kernel in $\R^{n+1}=\R^n_x\times\R_t$ and
set $(v_1,u_1,q_1)\equiv 0$ in $\R^n_x\times\R_t$.

\medskip

\noindent Having obtained $z_j:=(v_j,u_j,q_j)$ for $j\leq k$ and $\eta_1,\dots,\eta_{k-1}$ we choose 
\begin{equation}\label{eCI:eta}
\eta_k<2^{-k}
\end{equation}
in such a way that
\begin{equation}\label{eCI:osc}
\|z_k-z_k*\rho_{\eta_k}\|_{L^2(\Omega)}<2^{-k}.
\end{equation}
Then we apply Lemma \ref{l:onestep} to obtain $z_{k+1}=(v_{k+1},u_{k+1},q_{k+1})\in X_0$ such that
\begin{equation}\label{eCI:norm}
\|v_{k+1}\|^2_{L^2(\Omega)}\geq \|v_k\|^2_{L^2(\Omega)}+
\beta\left(|\Omega|-\|v_k\|^2_{L^2(\Omega)}\right)^2,
\end{equation}
and
\begin{equation}\label{eCI:weak}
\|(z_{k+1}-z_k)*\rho_{\eta_j}\|_{L^2(\Omega)}<2^{-k}\quad\textrm{ for all }j\leq k.
\end{equation}

\medskip

\noindent The sequence $\{z_k\}$ is bounded in $L^{\infty}(\R^n_x\times\R_t)$, therefore by passing to a suitable subsequence we may assume without loss of generality that
$$
z_k\overset{*}{\rightharpoonup} z\quad\textrm{ in }L^{\infty}(\R^n_x\times\R_t)
$$
for some $z=(v,u,q)\in X$, and that the sequence $\{z_k\}$ and the corresponding sequence $\{\eta_k\}$ satisfies the properties \eqref{eCI:eta},\eqref{eCI:osc},\eqref{eCI:norm} and \eqref{eCI:weak}. Then for every $k\in\N$
\begin{equation*}
\begin{split}
\|z_k*\rho_{\eta_k}-z*\rho_{\eta_k}\|_{L^2(\Omega)}
&\leq \sum_{j=0}^{\infty}\|z_{k+j}*\rho_{\eta_k}-z_{k+j+1}*\rho_{\eta_k}\|_{L^2(\Omega)}\\
&\leq \sum_{j=0}^{\infty}2^{-(k+j)}\leq 2^{-k+1},
\end{split}
\end{equation*}
and  since
\begin{equation*}
\begin{split}
\|z_k-z\|_{L^2(\Omega)}\leq 
&\|z_k-z_k*\rho_{\eta_k}\|_{L^2(\Omega)}+
\|z_k*\rho_{\eta_k}-z*\rho_{\eta_k}\|_{L^2(\Omega)}\\
&+\|z*\rho_{\eta_k}-z\|_{L^2(\Omega)},
\end{split}
\end{equation*}
we deduce that $v_k\to v$ strongly in $L^2(\Omega)$. 

Therefore, passing into the limit in \eqref{eCI:norm}
we conclude
\begin{equation}\label{e:final}
\|v\|^2_{L^2(\Omega)}\geq \|v\|^2_{L^2(\Omega)}+
\beta\left(|\Omega|-\|v\|^2_{L^2(\Omega)}\right)^2
\end{equation}
and hence $\|v\|^2_2 = |\Omega|$. Since $v$ vanishes
outside $\Omega$ and $|v|\leq 1$ in $\Omega$,
we conclude that $|v|= {\bf 1}_\Omega$.
Since $(v,u,q)\in X$, we also have that
$(v,u)(x,t)\in K^{co}$ for a.e. $(x,t)\in \Omega$. 
From this we deduce that 
$(v,u)(x,t)\in K$ for a.e. $(x,t)\in \Omega$, 
thus concluding the proof.
\end{proof}

\bibliographystyle{acm}

%\bibliography{eulerbib}

\end{document}